%
% slightly revised Sept. 2013
%   fixed error in proof of Lemma 2.2
%   and simple typos
%
\input amstex
\UseAMSsymbols
\NoBlackBoxes
\documentstyle{amsppt}
\magnification=\magstep1
\vsize=7.5in
\topmatter
\title{Chebyshev's conjecture and the prime number race}
\endtitle
\author Kevin Ford$^*$, Sergei Konyagin$^\dag$
\endauthor
\thanks
$^*$ First author supported in part by National Science Foundation grant
DMS-0070618.
\endthanks
\thanks
$^\dag$ Second author supported by INTAS grant N.~99-01080.
\endthanks
\rightheadtext{Prime number race and zeros of $L$-functions}
\leftheadtext{K. Ford, S. Konyagin}
\date July, 2002 \enddate
\noindent\address K.F: Department of Mathematics,
 University of Illinois at Urbana-Champaign,
 Urbana, IL 61801; \endaddress
\address
S.K.: Department of Mechanics and Mathematics, Moscow State University,
Moscow 119899, Russia.
\endaddress
\endtopmatter
\document
\predefine\barunder{\b}

\def\li{\text{\rm li}}
\def\sg{\sigma}
\def\g{\gamma}

\redefine\b{\beta}

\define\bc{\overline{\chi}}
\def\vp{\varphi}
\def\Li{\text{\rom{Li}}}
\def\curly{\Cal}
\def\BB{\curly B}

\redefine\le{\leqslant}
\redefine\ge{\geqslant}

\def\a{\alpha}
\def\e{\varepsilon}

\define\({\left(}
\define\){\right)}
\define\pfrac#1#2{\( \frac{#1}{#2} \)}
\define\bfrac#1#2{\left[ \frac{#1}{#2} \right]}
\define\ds{\displaystyle}

\head 1. Introduction and statement of problems \endhead

Dirichlet in 1837 proved that for any $a,q$ with $(a,q)=1$
there are infinitely many primes $p$ with $p\equiv a\pmod q$;
see, for example, ([Da], chapter 4). Chebyshev [Ch] noted
in 1853 that there are ``more'' primes congruent to $3$ than $1$
modulo $4$. His conjecture states that
$$\lim_{x\to\infty}\sum_{p>2}(-1)^{(p-1)/2}e^{-p/x}=-\infty.$$
As it was shown by Hardy, Littlewood [HL] and Landau [La2]
this holds if and only if the function
$$
L(s,\chi_1)=\sum_{n=0}^\infty\frac{(-1)^n}{(2n+1)^s}
$$
does not vanish for $\Re s>1/2$.

Chebyshev's conjecture was the origin for a big branch
of modern Number Theory, namely, comparative
prime--number theory. As usual, $p$ runs over the primes,
$k$ is a positive integer, $(k,l)=1$,
$$
\pi(x,k,l)=\sum\Sb p\le x\\p\equiv l\pmod k\endSb 1,
$$
$\pi(x)=\pi(x,1,1)$.  The functions
$\pi(x,k,l)$ for fixed $k$ and $(l,k)=1$ are all asymptotically
$x/(\vp(k)\log x)$ [Da, Ch. 20], where $\vp(k)$ is Euler's function, i.e.
the number of positive integers $l\le k$ with $(k,l)=1$. 
Still, as Chebyshev observed, there are
interesting inequities in the functions $\pi(x,k,l)$ for fixed $k$.
Of particular interest is the behavior of the functions
$$
\Delta(x,k,l_1,l_2)=\pi(x,k,l_1)-\pi(x,k,l_2).
$$
In Chebyshev's case, for example, $\Delta(x,4,3,1)$ is negative for the
first time at $x=26861$ [Lee].  More dramatically,
$\Delta(x,3,2,1)$ is negative for the first time at
$x=608,981,813,029$ [BH1].  Nonetheless,
Littlewood [Li] proved in 1914 that the functions $\Delta(x,4,3,1)$ and
$\Delta(x,3,2,1)$ each change sign infinitely often.
In their fundamental series of papers on comparative prime--number
theory ([KT1], [KT2]),
Knapowski and Tur\'an generalized Littlewood's theorem and
also indicated a large number of problems related to
comparison of $\pi(x,k,l_1)$ and $\pi(x,k,l_2)$.
Below we list some of these (and other) problems.
In the sequel discussing sign changes of $\Delta(x,k,l_1,l_2)$
and similar differences we shall assume that $l_1\not\equiv l_2\pmod k$
and $(l_1,k)=(l_2,k)=1$.

1. ``Infinity of sign changes''.
 To prove that $\Delta(x,k,l_1,l_2)$ changes sign infinitely often.

2.  ``Big sign changes''. To prove that
$\Delta(x,k,l_1,l_2)$ is $>x^{1/2-\e}$ (respectively $< - x^{1/2-\e}$)
for an unbounded sequence of $x$.
%arbitrary small $\e>0$ there is a sequence $x_1<x_2<\dots\to+\infty$
%such that, for $\nu=1,2,\dots,$
%$$\Delta(x_\nu,k,l_1,l_2)>x_\nu^{1/2-\e};$$
%and hence owing to the symmetry of $l_1,l_2$ also a sequence
%$y_1<y_2<\dots\to\infty$ such that
%$$\Delta(y_\nu,k,l_1,l_2)<-y_\nu^{1/2-\e}.$$
The use of function $x^{1/2-\e}$ is motivated by the fact
that if Extended Riemann Hypothesis for $k$ (see below)
is true than the inequality
$$
|\Delta(x,k,l_1,l_2)|=O(x^{1/2}\log x)\quad(x\ge 2)\tag{1.1}
$$
holds (see e.g. [Da, Ch. 20, (14)]).

3. ``Localized sign changes''. To prove that (i) for $T>T_0(k)$ and suitable
$G(T)<T$ the function $\Delta(x,k,l_1,l_2)$ changes sign
in the interval $G(T)\le x<T$; (ii-a) prove that
$\Delta(x,k,l_1,l_2)$ takes ``large'' positive and negative values
in the interval  $G(T)\le x<T$, where ``large'' means $> T^{1/2-\e}$;
(ii-b) find lower bounds on the number of sign changes for $x\le T$.

%4. ``Localized big sign changes''. To prove that for $T>T_1(k)$ and suitable
%$A(T)<T$ the inequalities
%$$\max_{A(T)\le x<T}\Delta(x,k,l_1,l_2)>\frac{T^{1/2}}{\Phi(T)}$$
%and hence also
%$$\min_{A(T)\le x<T}\Delta(x,k,l_1,l_2)<-\frac{T^{1/2}}{\Phi(T)}$$
%hold, with a $\Phi(x)>0$, satisfying also
%$$\lim_{x\to+\infty}\frac{\log\Phi(x)}{\log x}=0.\tag2$$

4. ``First sign change''. To determine a function $B(k)$ such that
for $1\le x\le B(k)$ all $\Delta(x,k,l_1,l_2)$ change sign at least once.

%6. ``Asymptotic estimation of the number of sign changes''. Clear.

5. ``Average preponderance problems''. A typical question is the following.
Denote by $N(x)$ the number of integers $n\le x$ with the property
$\Delta(x,4,3,1)\le 0$. Does the relation
$$
\lim_{x\to+\infty}\frac{N(x)}x=0\tag{1.2}
$$
hold?

%8. ``Strongly localized accumulation problems''. Is it true that
%for $T>C(k)$, for any $l_1,l_2$ with $l_1\not\equiv l_2\pmod k$
%and for suitable $T\le U_1<U_2\le2T$ we have
%$$
%\sum\Sb U_1\le p\le U_2\\ p\equiv l_1\pmod k\endSb 1
%-\sum\Sb U_1\le p\le U_2\\ p\equiv l_1\pmod k\endSb 1
%>\frac{\sqrt T}{\Phi(T)}
%$$
%where a function $\Phi(x)$ satisfies (2)?

6. ``Littlewood--generalizations''. A typical problem of this type would be
the existence of a sequence $x_1<x_2<\dots\to+\infty$ such that
simultaneously the inequalities
$$
\pi(x,4,1)\ge\frac12 \Li (x_\nu)
$$
and
$$
\pi(x,4,3)\ge\frac12 \Li (x_\nu)
$$
hold, where
$$
\Li(x) = \lim_{\e\to 0^+} \int_{0}^{1-\e} \frac{dt}{\log t}
 + \int_{1+\e}^x \frac{dt}{\log t}.
$$
This would constitute an obvious generalization of Littlewood's
theorem [Li] that for a suitable sequence $y_1<y_2<\dots\to+\infty$
the inequality $\pi(y_\nu)\ge \Li(y_\nu)$ holds.

7. ``Racing problems''.
%If $l_1,l_2,\dots,l_{\vp(k)}$ is any prescribed order of
%the reduced residue system $\pmod k$ then for a suitable
%sequence $x_1<x_2<\dots\to+\infty$ the inequalities
%$$
%\pi(x,k_\nu,l_1) > \pi(x,k_\nu,l_2) > \dots > \pi(x,k_\nu,l_{\vp(k)})
%$$
%hold.
The ``prime number race'' is colorfully described in [KT2, I]
in the following way. Consider a game with $r$ players, $2\le r\le \vp(k)$,
called ``$l_1$'' through ``$l_r$'' (here $l_1,\ldots l_r$ are
mutually incongruent modulo $k$), and at time $t$, each player ``$j$'' has
a score of $\pi(t,k,j)$ (i.e. player ``$j$'' scores 1 point whenever
$t$ reaches a prime $\equiv j \pmod{k}$).  As $t\to \infty$, will
each player take the lead infinitely often?  More generally,
will all $r!$ orderings of the players occur for infinitely many integers $t$?
This type of question originated in a paper of Shanks [Sh], who calculated
that $\pi(x,8,1)\le \max_{a\in \{3,5,7\}} \pi(x,8,a)$ for $x\le 10^6$.
It is generally believed that the answers to both questions are yes.
 If $r=2$ this is just the first problem
(infinity of sign changes).

8. ``Distribution problems''.  
Investigate the distribution of $\Delta(x,k,l_1,l_2)$, e.g.
study  $S(x;z)=\{1\le y\le x: 
\Delta(y,k,l_1,l_2)\le z (x^{1/2}/\log x)\}$ in a more general
way than indicated by problems 1--6.
For the general race problem, study distribution properties of the vectors
$(\pi(x,k,l_1),\ldots,\pi(x,k,l_r))$.

9. ``Union--problems''. For a given modulus $k$ and
disjoint subsets $A$ and $B$ of reduced residue classes mod $k$,
study the distribution of the function
$$
\sum_{p\in A, p\le x}1 - \frac{|A|}{|B|} \sum_{p\in B, p\le x}1. \tag{1.3}
$$
An important example is a generalization of Chebyshev's example:
let $A$ be the set of quadratic non-residues modulo $k$, and take $B$
to be the set of quadratic residues.

\bigskip

Besides the functions $\pi(x,k,l)$,
the distribution of primes in arithmetic progressions can be characterized
by some other functions which are more convenient to work
with. Let $\Lambda(n)$
be the Dirichlet--von Mangoldt function, namely, $\Lambda(n)=\log p$
if $n=p^m$ for some prime $p$ and some positive integer $m$, and $\Lambda(n)=0$
otherwise. The following functions are studied:
$$
\align
\psi(x,k,l)  &=\sum\Sb n\le x\\n\equiv l\pmod k\endSb \Lambda(n), \\
\Pi(x,k,l)   &=\sum\Sb n\le x\\n\equiv l\pmod k\endSb \frac{\Lambda(n)}
{\log n}, \\
\theta(x,k,l)&=\sum\Sb p\le x\\p\equiv l\pmod k\endSb \log p.
\endalign
$$
All problems 1--9 can be formulated for comparison of these functions.
Moreover, one can consider all questions ``in the Abelian sense''
as in original Chebyshev's paper [Ch]; for example, $\pi(x,k,l)$
can be replaced by $\sum_{p\equiv l\pmod k}e^{-pr}$ for $r>0$.
In this paper, however, we will be concerned only with problems involving the
functions $\pi(x,k,l)$, concentrating on problems 1 and 7.

%%%%%%%%%%%%%%%%%%%%%%%%%%%%%%%%%%%%%%%%%%%%%%%%%%%%%%%%%%%%%%%%%
%
\head 2. Analytic Tools \endhead
%
%%%%%%%%%%%%%%%%%%%%%%%%%%%%%%%%%%%%%%%%%%%%%%%%%%%%%%%%%%%%%%%%%

The methods of investigation of oscillatory properties
of the functions $\pi(x,k,l_1)-\pi(x,k,l_2)$ are similar to
those used to study oscillatory properties of the remainder term
of the prime number theorem, i.e. $E(x) = \pi(x)-\Li(x)$.
The primary tools are so-called ``explicit formulas'' linking
 the functions $\pi(x,k,l)$ 
to the distribution of the zeros in the critical
strip $0<\Re s < 1$ of the Dirichlet $L$-functions $L(s,\chi)$
for the characters $\chi$ modulo $k$. Note that
$$
\align
\sum_{\Sb n\le x \\ n\equiv l\pmod{k} \endSb} \frac{\Lambda(n)}{\log n}
&=\sum_{\Sb p\le x \\ p\equiv l\pmod{k} \endSb}
\frac{\Lambda(p)}{\log p}
+\sum_{\Sb p^2\le x \\ p^2\equiv l\pmod{k} \endSb}
\frac{\Lambda(p^2)}{\log (p^2)}+O(x^{1/3}) \\
&=\pi(x,k,l)+\sum_{\Sb 1\le u\le k \\ u^2\equiv l\pmod{k} \endSb}
\frac{\pi(x^{1/2},k,u)}2+O(x^{1/3}).
\endalign
$$
Using the asymptotic formula for $\pi(x,k,l)$ we have
$$
\pi(x,k,l)=\sum_{\Sb n\le x \\ n\equiv l\pmod{k} \endSb}
 \frac{\Lambda(n)}{\log n}
-\frac{N_k(l)}{\vp(k)}\frac{x^{1/2}}{\log x}
+o\left(\frac{x^{1/2}}{\log x}\right)\quad (x\to\infty)\tag{2.1}
$$
where $N_k(l)$ is the number of incongruent solutions of
the congruence $u^2\equiv l\pmod k$.

Let $D_k$ ($C_k$) denote the set of all (correspondingly, non-principal)
characters modulo $k$. For $\chi\in D_k$ define
$$\Psi(x;\chi) = \sum_{n\le x} \Lambda(n) \chi(n).$$
It is not difficult to show that
$$
\split
\vp(k)\sum_{\Sb n\le x \\ n\equiv l\pmod{k} \endSb} \frac{\Lambda(n)}{\log n}
&= \vp(k) \( \frac{\psi(x,k,l)}{\log x}+
  \int_2^x \frac{\psi(t,k,l)}{t\log^2 t}\, dt \) \\
&=\sum_{\chi\in D_k}\overline\chi(l)\left(\frac{\Psi(x;\chi)}{\log x}
  +\int_2^x\frac{\Psi(t;\chi)}{t \log^2 t}dt\right).
\endsplit
$$
Writing $D(x,k,l_1,l_2)=\psi(x,k,l_1)-\psi(x,k,l_2)$, we have by (2.1)
$$
\aligned
\vp(k)\Delta(x,k,l_1,l_2) &=\vp(k) \( \frac{D(x,k,l_1,l_2)}{\log x}
  +\int_2^x \frac{D(t,k,l_1,l_2)}{t\log^2 t}\, dt \) \\
&\qquad -(N_k(l_1)-N_k(l_2))\frac{x^{1/2}}{\log x}
  +o\left(\frac{x^{1/2}}{\log x}\right) \\
&=\sum_{\chi\in C_k}(\overline\chi(l_1)
  -\overline\chi(l_2))\left(\frac{\Psi(x;\chi)}{\log x}
  +\int_2^x\frac{\Psi(t;\chi)}{t \log^2 t}dt\right)\\
&\qquad -(N_k(l_1)-N_k(l_2))\frac{x^{1/2}}{\log x}
  +o\left(\frac{x^{1/2}}{\log x}\right)\quad (x\to\infty).
\endaligned\tag{2.2}
$$
By well-known explicit formulas (Ch. 19, (7)--(8) in [Da]),
when $\chi\in C_k$,
$$
\Psi(x;\chi) = - \sum_{|\Im \rho| \le x} \frac{x^\rho}{\rho} + O\(\log^2 x \),
\quad(x\ge2)\tag{2.3}
$$
where the sum is over zeros $\rho$ of $L(s,\chi)$ with $0<\Re \rho < 1$.

The zeros of $L(x,\chi)$ with largest real part will dominate the
sum in (2.3). 
The Extended Riemann Hypothesis for $k$ (abbreviated ERH$_k$)
states that all these zeros for all $\chi\in D_k$ lie on the critical line
$\Re s = \frac12$. Since
$$
\sum_{0\le \Im \rho \le T}1=O(T\log T)\quad(T\ge2)\tag{2.4}
$$
 for each $\chi$ ([Da], Ch. 16, (1)), one can see that
ERH$_k$ plus (2.3) implies (1.1). 
Moreover, it is known that if ERH$_k$ holds then
$$
\frac1x\int_0^x|\Psi(t;\chi)|^2 dt=O(x)
$$
for any $\chi\in C_k$ (this follows easily by the methods
in Chapter III of Cram\'er [Cr2]).  In this case $\Delta(x,k,l_1,l_2)$ has
average order at most $x^{1/2}/\log x$. If also
$N_k(l_1)\neq N_k(l_2)$ then the term
$$
-(N_k(l_1)-N_k(l_2))\frac{x^{1/2}}{\log x} \tag{2.5}
$$
in (2.2) is very significant, causing a shift in the mean value of
$\Delta(x,k,l_1,l_2)$.
In particular, we can expect that
if $l_1$ is a quadratic residue and $l_2$ is a quadratic nonresidue
modulo $k$ then for ``most'' $x$'s $\pi(x,k,l_1)<\pi(x,k,l_2)$.
Sometimes this phenomenon is called ``Chebyshev's Bias''.

A basic tool for proving oscillation theorems is the following
result of Landau [La1] on the location of singularities of the
Mellin transform of a non-negative function.

\proclaim{Lemma 2.1 [La1, p. 548]} 
Suppose $f(x)$ is real valued, and also non-negative for $x\ge x_0$.
Suppose also for some real numbers $\beta<\sigma$ that the Mellin transform
$$
g(s) = \int_1^{\infty} f(x) x^{-s-1}\, dx
$$
is analytic for $\Re s > \sigma$ and can be analytically continued to the
real segment $(\beta,\sigma]$.  Then $g(s)$ in fact represents an analytic
function in the half-plane $\Re s > \beta$.
\endproclaim

For example, if $f(x)=\vp(k) D(x,k,l_1,l_2)$ then
$$
g(s) = g(s,k,l_1,l_2) = -\frac{1}{s} \sum_{\chi\in C_k}
\(\bc(l_1)-\bc(l_2)\) \frac{L'(s,\chi)}{L(s,\chi)} \tag{2.6}
$$
for $\Re s>1$, with the right side providing a meromorphic continuation
of $g(s)$ to the whole complex plane.
Note that the poles of $g(s)$ (except $s=0$)
are a subset of the zeros of the functions
$L(s,\chi)$.  Also, $g(s)$ always has an infinite number of poles in
the critical strip [G1].
Assume that $g(s)$ has no real poles $s>\frac12$, $g(s)$ has a pole
$s_0$ with $\Re s_0 > \frac12$.  Take any $\a$ satisfying $\frac12 < \a
< \Re s_0$ and put $f(x)= (-1)^n D(x,k,l_1,l_2) + cx^\a$ for some 
constant $c$ and $n\in\{0,1\}$.
Applying Lemma 2.1 with different $n$ and $c$ we conclude that
$$
\limsup_{x\to \infty} \frac{D(x,k,l_1,l_2)}{x^\a} = \infty, \qquad
\liminf_{x\to \infty} \frac{D(x,k,l_1,l_2)}{x^\a} = -\infty.
$$
Since $\a > \frac12$, this is enough to deduce that $\Delta(x,k,l_1,l_2)$
changes sign infinitely often.  More generally, the following lemma
shows how oscillations of $D(x,k,l_1,l_2)$ and $\Delta(x,k,l_1,l_2)$
are related.

\proclaim{Lemma 2.2} Let
$$
h(x) = h(x,k,l_1,l_2) = \vp(k) x^{-1/2} D(x,k,l_1,l_2) - N_k(l_1)+N_k(l_2).
\tag{2.7}
$$
If
$$
\liminf_{x\to\infty} h(x) < 0 < \limsup_{x\to\infty} h(x),
$$
then $\Delta(x,k,l_1,l_2)$ has infinitely many sign changes.
\endproclaim

\demo{Proof}
By (2.2),
$$
\vp(k)\Delta(x,k,l_1,l_2) = \frac{x^{1/2}}{\log x}  h(x) + \int_2^x 
\frac{h(t)}{\sqrt{t}\, \log^2 t}\, dt +o\pfrac{\sqrt{x}}{\log x}
. \tag{2.8}
$$
By hypothesis, for an unbounded set of $x$ we have
$$
h(x) \ge \frac12 \max_{\log x\le y\le x} h(y).
$$
For such $x$ sufficiently large, (2.8) and the trivial bound
 $h(t)\ll \sqrt{t}$ (which we use for $t\le \log x$) imply that
$$
\vp(k)\Delta(x,k,l_1,l_2) \ge \frac{x^{1/2}}{\log x} h(x) \( 1 +O\( \frac{1}
{\log x} \)\) > 0.
$$
Similarly, there is an unbounded set of $x$ with $\Delta(x,k,l_1,l_2)<0$.
\enddemo

%%%%%%%%%%%%%%%%%%%%%%%%%%%%%%%%%%%%%%%%%%%%%%%%%%%%%%%%%%%%%%%%%
%
\head 3.  Sign changes in $\Delta(x,k,l_1,l_2)$ \endhead
%
%%%%%%%%%%%%%%%%%%%%%%%%%%%%%%%%%%%%%%%%%%%%%%%%%%%%%%%%%%%%%%%%%

\subhead 3.1.  The effect or large real zeros \endsubhead

By the remarks at the end of section 2,
to settle question 1 we need to deal with two cases: (i) $g(s)$
has a real pole $s>\frac12$; (ii) $g(s)$ has all of its poles in
the critical strip on the line $\Re s = \frac12$. 

Assume for the moment the following situation. Let $\chi_1\in C_k$ be
a real character, $L(\beta,\chi_1)=0$ for some $\beta\in(1/2,1)$,
but all other zeros of $L(s,\chi_1)$ and all zeros of $L(s,\chi)$,
$\chi\in D_k\setminus\{\chi_1\}$, have real part $\le \beta-\delta$
for some $\delta>0$.
Take $l$ such that $\chi_1(l)=-1$. It follows from (2.2)---(2.4)
that
$$
\vp(k)\Delta(x,k,1,l)=-\frac{2x^\beta}{\beta\log x}
+O\( \frac{x^{\beta}}{\log^2 x} \).
$$
Therefore $\pi(x,k,1)<\pi(x,k,l)$ for sufficiently large $x$.
This simple example shows that to prove the infinity of sign changes
of $\Delta(x,k,l_1,l_2)$ (and to succeed in any of the
problems 1--9)
for a general modulus $k$ we need some information
about the location of the zeros $\rho$.
Nowadays we cannot prove such properties in general, and they are
usually stated as suppositions.  Thus, many results of comparative
prime-number theory are conditional.  The most common supposition is
Haselgrove's condition: no $L(s,\chi)$, $\chi\in C_k$, vanishes
in the real interval $(0,1)$ (which eliminates the possibility of
case (i)).
To get effective results we need a bound on the
distance from any zero to the real axis,
and Haselgrove's condition is usually formulated as
the existence of an $A(k)$ such that no $L(s,\chi)$, $\chi\in C_k$,
vanishes in the parallelogram $0<\Re\rho<1$, $|\Im\rho|\le A(k)$.

\subhead 3.2. Infinitely many sign changes on Haselgrove's condition
\endsubhead

As already mentioned,
Littlewood [Li] proved that both functions
$\Delta(x,4,3,1)$ and
$\Delta(x,3,2,1)$ change sign infinitely often.
Knapowski and Tur\'an [KT1, II] extended these results significantly,
on the assumption of Haselgrove's condition.

\proclaim{Theorem 3.1 [KT1, II, Theorem 5.1]} If Haselgrove's condition
 is true for $k$ then the difference
$\pi(x,k,l)-\pi(x,k,1)$ changes sign infinitely often for any
$l$.
\endproclaim
In addition, Knapowski and Tur\'an
prove bounds for the first sign change as well as
for the magnitude and frequency
of the oscillations ([KT1, II, Theorems 5.1, 5.2], [KT1,III, Theorems
1.2, 1.3, 2.1, 3.1--3.4]).
In principal it is not difficult to verify Haselgrove's
condition for a particular $k$, and this has been done for
many small $k$ (P.C. Haselgrove, unpublished) including all $k\le 72$ [Ru1].
Recently, Conrey and Soundararajan [CS] proved that $L(s,\chi)\ne 0$
on $s\in (0,1)$
for at least 20\% of the real quadratic characters $\chi_{-8d}(n) = 
(\frac{-8d}{n})$, where $d$ runs over odd squarefree positive integers.
In [KT1, VI], a more general theorem is proved under
a slightly stronger hypothesis.  For some effectively
computable constant $c_1$ (independent of $k$), if no $L(s,\chi)$,
$\chi\in C_k$, vanishes in the domain $\Re s>1/2,\ |\Im s|\le c_1k^{10}$,
and $N_k(l_1)=N_k(l_2)$,
then $\Delta(x,k,l_1,l_2)$ changes sign infinitely often.
The last condition means that $l_1,l_2$
are simultaneously either quadratic residues or quadratic
nonresidues modulo $k$.
Later K\'atai [Ka1] proved the same conclusion under weaker hypotheses.

\proclaim{Theorem 3.2 [Ka1, Satz 2]}
Assuming Haselgrove's condition for $k$ and $N_k(l_1)=N_k(l_2)$,
then $\Delta(x,k,l_1,l_2)$ changes sign infinitely often.
More specifically, 
$$
\limsup_{x\to \infty} \frac{\Delta(x,k,l_1,l_2)}{\sqrt{x}/\log x} > 0,
\qquad \liminf_{x\to \infty} \frac{\Delta(x,k,l_1,l_2)}{\sqrt{x}/\log x} < 0.
$$
\endproclaim

For the general case of arbitrary $l_1$ and $l_2$, Haselgrove's condition
implies infinitely many sign changes of the difference
$D(x,k,l_1,l_2)$ [KT1, VII] (see also [Ka1, Satz 1]),
but the proof breaks down for
$\Delta(x,k,l_1,l_2)$, even if one assumes the full ERH$_k$.  Basically,
the proven magnitude of oscillations
of $\psi(x,k,l_1)-\psi(x,k,l_2)$ (of order $\sqrt{x}$ in [Ka1, Satz 1])
 are insufficient to overcome the term (2.4) appearing in (2.2).
As mentioned at the end of section 2, the difficulty is when all of
the singularities of $g(s)$ (a linear combination of functions
$\log L(s,\chi)$ in our case) have real part $= \frac12$.

%%%%%%%%%%%%%%%%%%%%%%%%%%%%%%%%%
%
\subhead 3.3. Infinitely many sign changes when $N_k(l_1) \ne N_k(l_2)$, 
$l_1\ne 1$, $l_2\ne 1$ \endsubhead
%
%%%%%%%%%%%%%%%%%%%%%%%%%%%%%%%%%

Fix $k,l_1,l_2$ and define $g(s)$ as in (2.6) and $h(x)$ as in (2.7).
Suppose $g(s)$ has no poles with $\Re s > \frac12$.  Label the poles
in the right half plane with positive imaginary part as $\frac12+i\g_1$,
$\frac12+i\g_2, \dots$ and let $G=\{ \g_1, \g_2, \ldots \}$.
Put $\g_0=0$ and let 
$a_0$ be the residue of $g(s)$ at $s=\frac12$ ($a_0=0$ if $g(s)$ is
analytic at $s=\frac12$) and let $a_j$ be the residue of $g(s)$
at $s=\frac12+i\g_j$ (typically the numbers $a_j$ have order $1/\g_j$).
Define 
$$
A(u) = \vp(k) e^{-u/2} D(e^u,k,l_1,l_2) = h(e^u)+N_k(l_1)-N_k(l_2)
$$
and for each $T>0$ let
$$
A_T^*(u) = \sum_{|\g_j| < T} \(1 - \frac{|\g_j|}{T} \) a_j e^{i\g_j u}
= a_0 + 2\Re \sum_{0<\g_j<T} \(1 - \frac{\g_j}{T} \) a_j e^{i\g_j u}.
$$
By Theorem 1 of Ingham [I],
$$
\liminf_{u\to\infty} A(u) \le \liminf_{u\to\infty} A_T^*(u) \le 
\limsup_{u\to\infty} A_T^*(u)
\le \limsup_{u\to\infty} A(u). \tag{3.1}
$$
As a consequence, if the numbers $\{ \g_i : 0 < \g_i < T\}$ 
are linearly independent
over the rationals then Kronecker's Theorem implies
$$
\align
\limsup_{u\to\infty} A_T^*(u) &=  a_0 + 2\Re \sum_{0<\g_j<T}
 \(1 - \frac{\g_j}{T} \) |a_j|, \\
\liminf_{u\to\infty} A_T^*(u) &=  a_0 - 2\Re \sum_{0<\g_j<T}
 \(1 - \frac{\g_j}{T} \) |a_j|.
\endalign
$$
In particular, if all $\g_1, \g_2,\ldots$ are linearly independent over
the rationals, then the above sum over $\g_j$ tends to $\infty$ as
$T\to\infty$.  In this case, by Lemma 2.2 the function $\Delta(x,k,l_1,l_2)$
changes sign infinitely often.

The linear independence property was introduced earlier
by Wintner [W1], [W2, chapter 13], and in [RS] it is called the
Grand Simplicity Hypothesis (GSH$_k$): The set of all $\gamma \ge 0$
such that $L(\frac12+i\gamma,\chi)=0$ for $\chi\in C_k$,
are linearly independent over $\Bbb Q$.
Note that GSH$_k$ implies that all the zeros are simple and that
$L(\frac12,\chi)\neq0$ for all such $\chi$.

\proclaim{Theorem 3.3 [I]}  Assume ERH$_k$ and GSH$_k$.  Then every
function $\Delta(x,k,l_1,l_2)$ changes sign infinitely often.
\endproclaim

It is not possible at the moment to show that any set of $\g_i$ are
linearly independent.  One way to avoid this is to use a computer to
find values of $u$ so that $A_T^*(u)$ is small or large and use the
inequality
$$
\liminf_{u\to\infty} A(u) \le A_T^*(u) \le \limsup_{u\to\infty} A(u),
$$
which follows from (3.1) and the almost periodic properties of $A_T^*(u)$.

Stark [St] generalized Ingham's theorem [St, Theorem 1],
in particular showing that $\Delta(x,5,4,2)$
changes sign infinitely often, which neither  Knapowski and  Tur\'an
nor K\'atai could prove with their methods.  A limiting case of Stark's
method allows one to prove results without any knowledge of the
non-trivial zeros of $L(s,\chi)$ [St, Theorem 3].
Stark also proved a more general theorem where one compares $\pi(x,k_1,l_1)$
with $\pi(x,k_2,l_2)$ with $k_1 \ne k_2$ (this theorem can
also be proved from the methods of [Ka1]).

\proclaim{Theorem 3.4 [St, Theorem 2 (i)]}
Suppose $(l_1,k_1)=(l_2,k_2)=1$, $N_{k_1}(l_1)=N_{k_2}(l_2)$
 and no $L$-function
in $C_{k_1} \cup C_{k_2}$ has a real zero $> \frac12$.  Then
$$
\limsup_{x\to\infty} \frac{\vp(k_1) \pi(x,k_1,l_1)-\vp(k_2)\pi(x,k_2,l_2)}
{\sqrt{x}/\log x} > 0.
$$
\endproclaim

About the same time, Bateman et al [BBHKS] showed that the
linear independence condition in Ingham's theorem could be relaxed
considerably with only a slightly weaker conclusion.  Diamond [Di]
extended and generalized this result.  To state it, we introduce
a notion of weak independence.  Let $N$ be a positive integer.
A subset of $G$, $\{\g_{j_1}, \g_{j_2}, \ldots,
\g_{j_m} \}$, is said to be $N$-independent if for every choice of integers
$n_1, \ldots, n_m$ satisfying $|n_r| \le N$ for each $r$ and
$\sum_{r=1}^m |n_r| \ge 2$ we have $\sum_{r=1}^m n_r \g_{j_r} \not\in G$.
With modest values of $m$ and $N$, this can be checked by computer in
a reasonable time.  

\proclaim{Theorem 3.5 [Di]}
Suppose that  $\{\g_{j_1}, \g_{j_2}, \ldots,\g_{j_m} \}$ is $N$-independent.
Then
$$
\limsup_{u\to\infty} A(u) \ge a_0 + \frac{2N}{N+1} \sum_{r=1}^m |a_{j_r}|,
\quad \liminf_{u\to\infty} A(u) \le a_0 - \frac{2N}{N+1} \sum_{r=1}^m
|a_{j_r}|.
$$
\endproclaim

In light of Lemma 2.2, for every $k,l_1,l_2$ it is possible to
prove with a finite computation
that $\Delta(x,k,l_1,l_2)$ changes sign infinitely often
using Diamond's theorem.  Grosswald [G3] carried out such computations
for $k\in \{5, 7, 11, 13, 19\}$, in each case proving unconditionally
that all functions $\Delta(x,k,l_1,l_2)$ change sign infinitely often.
The cases $k=3,4$ and $6$ were settled by Littlewood [Li], $k=8$ by
Knapowski and Tur\'an [KT1], and earlier the cases $k\in \{43, 67, 163\}$
were settled by Grosswald [G2] (important here is the fact that
$a_1$ is quite large; this is connected to the class number of $\Bbb Q(
\sqrt{-k})$ being 1 [BFHR]).
We know of no computations for other $k$.  

Applying Theorem 3.5 with $N=1$ and $m=1$ gives the unconditional bounds
$$
\limsup_{u\to\infty} A(u) \ge a_0 + \max_{j\ge 1} |a_j|,
\quad \liminf_{u\to\infty} A(u) \le a_0 - \max_{j\ge 1} |a_j|.
$$

Kaczorowksi [K2] proposed different method for avoiding full
linear independence of the $\g_i$.  Let $\Omega=T^{\infty}$ denote
the infinite dimensional torus, i.e. the topological product of infinitely
many copies of $T=\{z\in \Bbb C: |z|=1\}$.  Define the continuous
homomorphism $\lambda(t)=(e^{2\pi i \g_1 t}, e^{2\pi i \g_2 t}, \ldots)$.
Let $\Gamma=\overline{\lambda(\Bbb R)}$, the closure of $\lambda(\Bbb R)$
in the Tikhonov topology of $\Omega$.  For $\a=(\a_1,\a_2,\ldots)\in 
\Omega$ define $u\a = (u\a_1,u\a_2,\ldots)$ and let $S_0$ be the
stabilizer of $\Gamma$, i.e. $S_0=\{u\in T: u\Gamma=\Gamma\}$.
As $S_0$ is a closed subgroup of $T$, either $S_0=T$ or $S_0$
is a cyclic group generated by a root of unity.  If the numbers
$\g_j$ are linearly independent, Kronecker's theorem implies
that $S_0=T$.  However, Kaczorowski proves that $h(x,k,l_1,l_2)$
takes arbitrarily large and small values on the hypothesis that $S_0$
contains any element other that $\pm 1$ ([K2], Corollaries 4,5).
This will be discussed further in \S 4.

%%%%%%%%%%%%%%%%%
%
\subhead 3.4.  First sign change \endsubhead
%
%%%%%%%%%%%%%%%%%%%%%%

The problem of bounding the first sign change of $\Delta(x,k,l_1,l_2)$
can be attacked in a brute-force way by computing all of the primes up
to a given limit using the sieve of Eratosthenes and then tabulate the
functions $\pi(x,k,l)$ directly.  Using modern computers, John Leech [Lee]
discovered in 1957 the first negative value in $\Delta(x,4,3,1)$, which
occurs at $x=26861$.  It turns out that when $k|24$, $(l,k)=1$ and
$1<l<k$ that negative values of $\Delta(x,k,l,1)$ are quite rare.
In a massive computation in the 1970s, Bays and Hudson ([BH1]--[BH4]) computed
these functions up to $x=10^{12}$.  In addition to several sign changes
in $\Delta(x,4,3,1)$, they discovered sign changes in only three other
such functions: $\Delta(x,3,2,1)$ with its first negative value at
$x=608,981,813,029$, $\Delta(x,8,5,1)$ with its first negative value at
$x=588,067,889$, and $\Delta(x,24,13,1)$ with its first negative
value near $978,000,000,000$.

Computations with the explicit formulas (2.3) but truncated to
a sum over $\rho$ with $|\Im \rho| \le T$ for a fixed $T$ (e.g. $T=10000$)
provides one with a good sense of where the sign changes likely
occur, since the smaller zeros of (2.3) contribute the dominant part
of the sum. 
Actually proving that a negative value occurs at or near $x$
can be done in three ways.  First, if $x$ is not too large (say $x<10^{14}$
with the latest computers), the above mentioned brute force method
can work.  For intermediate $x$ (up to maybe $10^{20}$ with today's
technology), one can use Hudson's [H1] generalization of the famous
Meissel formula to compute exactly the value of the functions
$\pi(x,k,l)$.  As shown by Lagarias, Miller and Odlyzko [LMO],
this can be done in time $O(x^{2/3+\e})$, compared to time $O(x^{1+\e})$
required for the brute-force method.  The first author wrote
a computer program implementing the algorithm, and one
result is that $\Delta(1.9282\times 10^{14},8,1,7)<0$.

The third method, which is the only practical one for really large
$x$, is based on the explicit formulas (2.3).
The first papers on the subject are those of Skewes ([Sk1], [Sk2]),
 who gave
in 1955 the first unconditional upper bound on the first sign change
of $\pi(x)-\Li(x)$, namely $10^{10^{10^{10^3}}}$.
Lehman [Leh] showed that the results of the computations of the 
greatly truncated sums  (2.3) for the Riemann zeta function
could be made rigorous, in the sense that  one could use the 
calculations to prove that $\pi(x)>\Li(x)$ somewhere
in a short interval of the form $[x_0,(1+\delta)x_0]$.  
Lehman's method drastically improved the upper bound for the first
sign change.  In particular, he proved that it must occur before
$1.5926\times 10^{1165}$  and using more known zeros for the Riemann
zeta function his method was used by te Riele [tR]
to lower the bound to $6.6658\times 10^{370}$ and by Bays and Hudson [BH5]
to lower it further to $1.39822\times 10^{316}$.
Ford and Hudson [FH] generalized Lehman's method, which allowed them to
localize sign changes of any linear combination of functions
$\pi(x,k,l)$, in particular the functions $\Delta(x,k,l_1,l_2)$.
A consequence is the following theorem, proved with the aid of
the zeros $\rho$ of $L$-functions with $|\Im \rho| \le 10000$ ([Ru1],[Ru2]).

\proclaim{Theorem 3.6 [FH, Corollary 4]}  
For each $b\in \{3,5,7\}$, $\pi_{8,b}(x) < \pi_{8,1}
(x)$ for some $x<5\times 10^{19}$.  For each $b\in \{5,7,11\}$,
$\pi_{12,b}(x)<\pi_{12,1}(x)$ for some $x<10^{84}$.
For each $b\in \{5,7,11,13,17,19,23\}$,
$\pi_{24,b}(x) < \pi_{24,1}(x)$ for some $x<10^{353}$.
\endproclaim

%%%%%%%%%%%%%%%%%%%%%%%%%%%%%%%%%%%%%%%%%%%%%%%%%%%%%%%%%%%%%%%%%%%%%%%%
%
\head 4. Racing problems  \endhead
%
%%%%%%%%%%%%%%%%%%%%%%%%%%%%%%%%%%%%%%%%%%%%%%%%%%%%%%%%%%%%%%%%%%%%%%%%

In racing problems for $r>2$ virtually nothing is known
unconditionally.  All results we know
are proven under ERH$_k$; sometimes other assumptions are used.
In the end of this paper we will explain principal difficulties to get
unconditional results in racing problems.

\subhead 4.1 Conditional Results \endsubhead

In [K3] Kaczorowski assuming ERH$_k$ solves the racing problem in a weak
form: he shows that $1\pmod k$ wins and loses infinitely often.

\proclaim{Theorem 4.1 [K3, Theorem 1]}
Let $k\ge 3$ and assume ERH$_k$. Then there exist
infinitely many integers $m$ with
$\ds \pi(m,k,1)>\max_{l\not\equiv1\pmod k}\pi(m,k,l)$.
 Moreover, the set of $m$'s
satisfying this inequality has positive lower density. The same statement
holds true for $m$ satisfying the inequality
$\ds \pi(m,k,1)<\min_{l\not\equiv1\pmod k}\pi(m,k,l)$.
\endproclaim

In [K5] Kaczorowski showed that certain other orderings of the functions
$\pi(x,k,l)$ occur for a set of $x$ having positive lower density, but
these results are very complicated and we do not state them here.
Note that Theorem 4.1 conditionally solves an average preponderance
problem for
$\pi(n,4,1)-\pi(n,4,3)$ as it stated in problem 5 of section 1, (1.2).
It shows that under ERH$_4$,
$$
\liminf_{x\to+\infty}\frac{N(x)}x>0.
$$ 
This inequality could also be deduced from Wintner [W1].

In [K3] Kaczorowski assuming ERH$_k$ solves the racing problem modulo $5$
for the function $\psi$. So, for any permutation $(l_1,l_2,l_3,l_4)$ of the
sequence $(1,2,3,4)$ there exist infinitely many integers $m$ with
$$
\psi(m,5,l_1)>\psi(m,5,l_2)>\psi(m,5,l_3)>\psi(m,5,l_4).\tag{4.1}
$$
No doubt, more extensive (maybe, technically not feasible)
calculations can give the same assertion for the function $\pi(x,5,l)$.
Moreover, one can try do the same job for other moduli.

All results on sign changes of the functions $\Delta(x,k,l_1,l_2)$ make
use of some kind of almost periodicity of the sums in (2.3).  For example,
the work of Knapowski and Tur\'an depends on results from Tur\'an's power
sum methods [T].
Kaczorowski takes a different approach, replacing $x$ in (2.3) by
$e^z$ where $z$ is a complex variable.
This gives the so-called $k$-functions [K1].
For $\Im z>0$ and $\chi\in D_k$ he defines
$$
k(z,\chi)=\sum_{\Im\rho>0}e^{\rho z}
$$
and
$$
K(z,\chi)=\int_{i\infty}^z k(s,\chi)ds=\sum_{\Im\rho>0}\frac1\rho e^{\rho z}
$$
where the summation is taken over all nontrivial zeros of $L(s,\chi)$
with positive imaginary parts.  Further, for $(l,k)=1$, $0<l\le k$
the following functions are defined
$$
F(z,k,l)=-2e^{-z/2}\frac1{\vp(k)}\sum_{\chi\pmod k}\overline\chi(l)
K(z,\chi')-\frac2{\vp(k)}\sum_{\chi\pmod k}\overline\chi(l)m(\frac12,\chi),
$$
where $\chi'$ denotes the primitive Dirichlet's character induced by
$\chi$ and $m(\frac12,\chi)$ is the multiplicity of a zero of $L(s,\chi)$
at $s=\frac12$ ($m(\frac12,\chi)=0$ when $L(\frac12,\chi)\neq0$).

Historically, the $k$-functions were introduced by Cram\'er ([Cr1], [Cr2])
in connection with the Riemann zeta function and the Dedekind zeta
functions.

For real $x$ let $P(x,k,l)=\lim_{y\to0+}\Re F(x+iy,k,l)$.  It can be proved
that the limit exists for all real $x$, and that $P(x,k,l)$ is
piecewise continuous.  Moreover, we have (cf. (2.3))
$$
P(x,k,l) = e^{-x/2} \( \psi(e^x,k,l)-\frac{e^x}{\phi(k)} \) 
+ E(x,k,l),
$$
where $|E(x,k,l)| = O(xe^{-x/2})$ ($x\ge 1$).
Even though $E(x,k,l)$ is very small, its behavior
near $x=0$ is important to many of the applications. 
Important in Kaczorowski's analysis is the property that for
fixed $y=\Im z >0$ each function $F(x+iy,k,l)$ is almost periodic in
the sense of Bohr, and further the functions
$P(x,k,l)$ are almost periodic in the sense of Stepanov (see e.g. [Be]
for definitions and properties of various types of
 almost periodic functions).
A consequence ([K3], Lemma 3) is that for every number $x_0$ and
every $\e>0$,
the vector $V(x)=(P(x,k,l_1),\ldots,P(x,k,l_n))$ satisfies
$\| V(x)-V(x_0) \| \le \e$ for a set of $x$ having positive lower density.
Thus, finding a single point $x_0$ with
$$
P(x_0,k,l_{j_1}) > P(x_0,k,l_{j_2}) > \cdots > P(x_0,k,l_{j_n})
$$
implies that the same inequality occurs for a set of $x$ having
positive lower density.  Computing $P(x,k,l)$ is easy for small $x$,
and to prove Theorem 4.1, Kaczorowski
finds appropriate $x_0$'s in a small neighborhood of $0$.

The proof of Kaczorowski's theorem mentioned in \S 3.3
uses the fact that the function $F(z,k,l)$ has logarithmic
singularities at the points $\pm \log(p^m)$ where $p^m\equiv l \pmod{k}$,
and thus $\Im F(z,k,l)$ is unbounded in a vicinity of these points.
Adopting the notation from \S 3.3, let $\zeta \in S_0$, $\zeta\ne \pm 1$.
Then, for fixed $y>0$ the function $F(x+iy,k,l)$ exhibits a kind of
``twisted'' almost periodicity, namely for every $\e>0$ there is a
number $w$ so that 
$$
\sup_{x\in \Bbb R} |F(x+iy,k,l)-\zeta F(x+w+iy,k,l)| \le \e.
$$
Since $\zeta\ne \pm 1$, this implies that the real part of $F(z,k,l)$
is unbounded in the upper half plane, and this can be used to prove
that  $P(x,k,l)$ is unbounded.

The paper [K6] contains a nice overview of Kaczorowski's methods and
results.  In addition, one can define $k$-functions for a wide class
of $L$-functions occurring in number theory and use them to study
oscillations and distribution of various ``error terms'' in prime
number theory [KR].

%Consider the case $k=5$. There is a formula
%$$
%P_1(x)=e^{-x/2}(\psi_0(e^x,5,1)+Q(x)+\frac12\log(e^x-1)+B_1).
%$$
%Here $B_1$ is a computable constant,
%$$\psi_0(x,k,a)=\frac12(\psi(x+0,k,a)+\psi(x-0,k,a)),$$
%$$Q(x)=\frac14\left(\log5\left[\frac x{\log5}\right]_0-x-e^x\right),$$
%where
%$$[x]_0=\frac12([x+0]+[x-0]).$$
%Similar expressions are written for $j=2,3,4$.
%On one hand,
%we can compute $P_j(x)$ for small $x$. On the other hand,
%we get asymptotic behavior of $\psi(e^x,5,j)$:
%$$e^{-x/2}(\psi(e^x,5,j)-e^x/4)=P_j(x)+o(1)\quad(x\to\infty).$$
%So, if we want to prove that the inequalities (2.4)
%can occur for an arbitrary large $m$, it is enough
%to check that for some positive $\delta>0$ there exists
%an arbitrary large $x$ with
%$$P_{l_j}(x)>P_{l_{j+1}}(x)+\delta\quad(j=1,2,3).\tag9$$
%To prove this, it is sufficient to find some $x$ satisfying (9) and
%to use that under ERH$_5$ the functions $P_j$ are Bohr almost periodic.
%For any permutation $(l_1,l_2,l_3,l_4)$ of the set $(1,2,3,4)$
%such $x$ is explicitly indicated.

The racing problem has been extensively investigated by
Rubinstein and Sarnak [RS]. Let $P\subset(0,+\infty)$,
$$
\split
\overline\delta(P) &= \limsup_{X\to\infty}\frac1{\log X}
\int_{t\in P\cap[2,X]} \frac{dt}t, \\
\underline\delta(P) &= \liminf_{X\to\infty}\frac1{\log X}
\int_{t\in P\cap [2,X]} \frac{dt}t.
\endsplit
$$
If the latter two quantities are equal, 
the logarithmic density $\delta(P)$ of the set $P$
is their common value.  The problem is to study the existence of
and to estimate the logarithmic density of the set
$$
P_{k;l_1,\dots,l_r}=\{x\ge2:\pi(x,k,l_1)>\pi(x,k,l_2)>\dots
>\pi(x,k,l_r\}.
$$
Introduce the vector valued functions
$$
E_{k;l_1,\dots,l_r}(x)=\frac{\log x}{\sqrt x}(\vp(k)\pi(x,k,l_1)-\pi(x),
\dots,\vp(k)\pi(x,k,l_r)-\pi(x)).
$$

\proclaim{Theorem 4.2 [RS, Theorem 1.1]} Assume ERH$_k$. 
Then $E_{k;l_1,\dots,l_r}$ has a limiting
distribution $\mu_{k;l_1,\dots,l_r}$ on $\bold R^r$, that is
$$
\lim_{X\to\infty}\frac1{\log X}\int_2^X f(E_{k;l_1,\dots,l_r}(x))\frac{dx}x
=\int_{\bold R^r}f(x)d\mu_{k;l_1,\dots,l_r}
$$
for all bounded continuous functions $f$ on $\bold R^r$.
\endproclaim

If it turns out that if the measure $\mu_{k;l_1,\dots,l_r}$ is
 absolutely continuous then
$$
\delta(P_{k;l_1,\dots,l_r})=\mu_{k;l_1,\dots,l_r}
(\{x\in\bold R^r:x_1>\dots>x_r\}).
$$
However, the authors of [RS] write that assuming only ERH$_k$ they
don't know that $\delta(P_{k;l_1,\dots,l_r})$ exists.

The measures $\mu$ are very localized but not compactly supported.
Let $B_R'=\{x\in\bold R^r:|x|>R\}$, $B_R^+=\{x\in\bold B_R':
\e(l_j)x_j>0\}$, $B_R^-=-B_R^+$, where $\e(l)=1$ if $l\equiv1\pmod k$
and $\e(l)=-1$ if $l\not\equiv1\pmod k$.

\proclaim{Theorem 4.3 [RS, Theorem 1.2]} 
Assume ERH$_k$. Then there are positive constants
$c_1,c_2,c_3,c_4$, depending only on $k$, such that
$$\mu_{k;l_1,\dots,l_r}(B_R')\le c_1\exp(-c_2\sqrt R),$$
$$\mu_{k;l_1,\dots,l_r}(B_R^{\pm})\ge c_3\exp(-\exp(c_4R)).$$
\endproclaim
The second inequality gives a quantitative version of
the theorem from [K3].

 Montgomery [Mo],
using ERH$_1$ (the Riemann Hypothesis) and GSH$_1$,
 investigated the tails of $\mu_1=\mu_{1:1}$.  He showed that
$$\exp(-c_2\sqrt R\exp\sqrt{2\pi R})\le\mu_1(B_R^{\pm})
\le\exp(-c_1\sqrt R\exp\sqrt{2\pi R}).$$
Rubinstein and Sarnak [RS] under ERH$_k$ and GSH$_k$
 have found an explicit formula
for the Fourier transform of $\mu_{k;l_1,\dots,l_r}$. Special cases of
the formula were proven earlier in [W1] and [Ho]. The formula implies
that, for $r<\vp(k)$, $\mu_{k;l_1,\dots,l_r}=f(x)dx$ with a rapidly decreasing
entire function $f$.  As a consequence, under ERH$_k$ and GSH$_k$ each
$\delta(P_{k;l_1,\dots,l_r})$ does exist and is nonzero (including the case
$r=\vp(k)$). Hence, conditionally the solution to the racing problem
is affirmative.

Also in [RS] is a procedure for calculating $\delta(P_{k;l_1,\ldots,l_r})$
using known zeros of $L$-functions ([Ru1],[Ru2]).  In particular,
they compute $\delta(P_{4;3,1}) = 0.9959\ldots$, thus giving a 
quantitative version of Chebyshev's statement.  Although for small $k$ many of
the densities are quite large,  as $k\to \infty$ they become more
uniform.

\proclaim{Theorem 4.4 [RS, Theorem 1.5]}  Assume ERH$_k$ and GSH$_k$
for all $k\ge 3$.  For $r$ fixed,
$$
\max_{l_1,\ldots,l_r} \left| \delta(P_{k;l_1,\ldots,l_r}) - \frac{1}{r!}
\right| \to 0 \text{ as } k\to \infty.
$$
\endproclaim

We say that
$k;l_1,\dots,l_r$ is unbiased if $\mu_{k;l_1,\dots,l_r}$ is
invariant under permutations of $(x_1,\dots,x_r)$. In this case
$\delta(P_{k;l_1,\dots,l_r})=1/r!$. 

\proclaim{Theorem 4.5 [RS, Theorem 1.4]} Assume ERH$_k$ and GSH. Then $k;l_1,\dots,l_r$
is unbiased if and only if either $r=2$ and $N_k(l_1)=N_k(l_2)$, or
$r=3$ and $l_2\equiv l_1g\pmod k$, $l_3\equiv l_1g^2\pmod k$, where
$g^3\equiv1\pmod k$.
\endproclaim

Feuerverger and Martin [FM] computed numerous densities for small
moduli $k$, studying in particular the cases where $r\ge 3$, $k;l_1,\ldots,
l_r$ is not unbiased and $N_k(l_1)=\cdots=N_k(l_r)$, e.g. $8;3,5,7$.
Curious inequities occur in these cases, which can be ``explained''
in terms of distribution properties of the functions $\Psi(x,\chi)$
(see (2.3) and also [Ma]).
Computation of the densities when $r\ge 3$ is very complex
using the methods in [RS] or [FM].  A much simpler and faster method
(but less rigorous) is given in [BFHR].

Rubinstein and Sarnak also studied the union problem mentioned in 
problem 9 in section 1.  Suppose $k\ge 3$, let $A$ be the set of quadratic
non-residues modulo $k$ and let $B$ be the set of quadratic residues.
Define
$$
P_{k;N,R} = \left\{ x\ge 2: \sum_{p\in A, p\le x} 1 > \frac{|A|}{|B|}
\sum_{p\in B, p\le x} 1 \right\}. 
$$
With $k$ restricted to integers possessing a primitive root
(again assuming ERH$_k$ and GSH$_k$), in [RS] it is proved that
$\delta(P_{k;N,R}) \to \frac12$ as $k\to\infty$,  although the
convergence is far from monotone.  In fact there is a connection 
between the values of $\delta(P_{k;N,R})$ and the class number of the
imaginary quadratic field $\Bbb Q(\sqrt{-d})$ [BFHR].

%%%%%%%%%%
%
\subhead 4.2 Barriers to unconditional results \endsubhead
%
%%%%%%%%%%%%%

One may ask if a racing problem for $r>2$ may be solved without the
assumption of ERH$_k$, since the problem of proving infinitely many sign
changes of $\Delta(x,k,l_1,l_2)$ is easier if
$g(s,k,l_1,l_2)$ (cf. end of section 2)
 has a pole with real part $>1/2$ and no real poles $>1/2$.  
Since it is believed that the sets of zeros for each $L(s,\chi)$ are
disjoint, this condition essentially says that ERH$_K$ is false.

In particular, can it be shown, for some quadruples $(k,l_1,l_2,l_3)$,
that the 6 orderings of the functions $\pi(x,k,l_j)$
occur for infinitely many integers $x$,
without the assumption of ERH$_k$ (while still
allowing the assumption of Haselgrove's condition and/or
that zeros with imaginary part $< R_k$
lie on the critical line for some constant $R_k$)? Recently we have
answered this question in the negative (in a sense) for all quadruples
$(k,l_1,l_2,l_3)$ [FK].  Thus, in a sense the hypothesis
ERH$_k$ is a necessary condition for proving
any such results when $r>2$.

Let $D=(k,l_1,l_2,l_3)$, where $l_1,l_2,l_3$
 are distinct residues modulo $k$ which are coprime to $k$.
Suppose for each $\chi\in C_k$, $B(\chi)$ is a sequence of complex numbers
with positive imaginary part (possibly empty, duplicates allowed), and
denote by $\BB$ the system of $B(\chi)$ for $\chi\in C_k$.  Let $n(\rho,\chi)$
be the number of occurrences of the number $\rho$ in $B(\chi)$.
The system $\BB$ is called a {\it barrier} for $D$ if the following hold:
(i) all numbers in each $B(\chi)$ have real part in $[\beta_2,\beta_3]$,
where $\frac12 < \beta_2 < \beta_3\le 1$; (ii)
for some $\beta_1$ satisfying $\frac12 \le \beta_1 < \beta_2$,
if we assume that
for each $\chi\in C_k$ and $\rho\in B(\chi)$, $L(s,\chi)$ has a zero
of multiplicity $n(\rho,\chi)$ at $s=\rho$, and all other zeros of $L(s,\chi)$
in the upper half plane have real part $\le \beta_1$, then one of the six
orderings of the
three functions $\pi(x,k,l_j)$ does not occur for large $x$.
If each sequence $B(\chi)$ is finite, we call $\BB$ a {\it finite
barrier} for $D$ and denote by $|\BB|$ the sum of the number of elements of each
sequence $B(\chi)$, counted according to multiplicity.

\proclaim{Theorem 4.6 [FK]}
For every real numbers $\tau>0$ and $\sigma>\frac12$ and
every $D=(k,l_1,l_2,l_3)$, there is a finite
barrier for $D$, where each sequence
$B(\chi)$ consists of numbers with real part $\le \sigma$ and imaginary
part $>\tau$.  In fact, for most $D$, there is a barrier with $|\BB|\le 3$.
\endproclaim

We do not claim that the falsity of ERH$_k$ implies that one of the six
orderings does not occur for large $x$.  For example, take $\sg > \frac12$,
and suppose each non-principal character modulo  $k$ has a unique zero
with positive imaginary part to the right of the critical line, at
$\sg + i \gamma_\chi$.  If the numbers $\gamma_\chi$ are linearly independent
over the rationals, it follows easily from Lemma 4.1 below and
Kronecker's Theorem that in fact
all $\phi(k)!$ orderings of the functions $\{ \pi(x,k,l) : (l,k)=1 \}$
occur for an unbounded set of $x$.

As an example,
we will demonstrate the existence of a barrier in a more simple situation:
$r=4$, $k=5$, $l_j=j\ (j=1,2,3,4)$.
By (2.2)--(2.4) we have the following lemma.

\proclaim{Lemma 4.1} Let $\beta\ge \frac12$, $x\ge 10$
and for each $\chi\in C_k$,
let $B(\chi)$ be the sequence of zeros (duplicates allowed)
of $L(s,\chi)$ with  $\Re s>\beta$ and $\Im s > 0$.   Suppose further that
all $L(s,\chi)$ are zero-free on the real segment $0<s<1$.  If
$(l_1,k)=(l_2,k)=1$ and $x$ is sufficiently large, then
$$
\phi(k)\Delta(x,k,l_1,l_2)= -2\Re \left[ \sum_{\chi\in C_k}
(\bc(l_1)-\bc(l_2)) \sum_{\Sb \rho\in B(\chi) \\ |\Im \rho|\le x \endSb}
 f(\rho) \right] + O(x^\b \log^2 x),
$$
where
$$
f(\rho) := \frac{x^{\rho}}{\rho\log x} + \frac{1}{\rho}
\int_2^x \frac{t^{\rho}}{t\log^2 t}
\, dt= \frac{x^{\rho}}{\rho\log x} + O\pfrac{x^{\Re \rho}}{|\rho|^2\log^2 x}.
$$
\endproclaim
Here the constant implied by the Landau $O$-symbol may depend on $k$,
but not on any other variable.
Take $\chi_1\in C_5$ so that $\chi_1(1)=1, \chi_1(2)=i, \chi_1(3)=-i,
\chi_1(4)=-1$. Let $t$ be a large positive number.
Take $\sigma>\frac12$, $B(\chi_1)=\{\sigma+it\}$,
$B(\chi)=\emptyset$ for $\chi\in C_5\setminus\{\chi_1\}$.
We use Lemma 4.1 with $\frac12 < \beta<\sigma$. For $\rho=\sigma+it$ we have
$$
f(\rho)=\frac{-x^\sigma}{t\log x}x^{it}i
+O\left(\frac{x^\sigma}{t^2\log x}\right).\tag{4.2}
$$
We claim that in this situation the inequality
$$
\pi(x,5,1)>\pi(x,5,4)>\pi(x,5,2)>\pi(x,5,3)\tag{4.3}
$$
cannot occur for large $x$. Indeed, by Lemma 4.1 and (4.2), for large $x$
(depending on $\sigma$, $\beta$) we have
$$
\pi(x,5,1) - \pi(x,5,4) = \frac{x^\sigma}{t\log x}\Re(x^{it}i)
+O\left(\frac{x^\sigma}{t^2\log x}\right),\tag{4.4}
$$
$$
\pi(x,5,4) - \pi(x,5,2) = \frac{x^\sigma}{t\log x}\Re(x^{it}(-1/2-i/2))
+O\left(\frac{x^\sigma}{t^2\log x}\right),\tag{4.5}
$$
$$
\pi(x,5,2) - \pi(x,5,3) = \frac{x^\sigma}{t\log x}\Re(x^{it})
+O\left(\frac{x^\sigma}{t^2\log x}\right).\tag{4.6}
$$
It is not difficult to prove that for any $t$
$$
\min(\Re(x^{it}i),\Re(x^{it}(-1/2-i/2)),\Re(x^{it}))\le
-\sqrt{0.1}.
$$
Therefore, the estimates (4.4)--(4.6) imply
$$
\gather
\min(\pi(x,5,1) - \pi(x,5,4),\pi(x,5,4) - \pi(x,5,2),
\pi(x,5,2) - \pi(x,5,3))\\
\le-\frac{x^\sigma}{\sqrt{10}t\log x}
+O\left(\frac{x^\sigma}{t^2\log x}\right).
\endgather
$$
For large $t$ this does not agree with (4.3).

For $r=3$ our construction of finite barriers is more difficult.
It uses multiple zeros or zeros of several functions $L(s,\chi)$
simultaneously. However, answering a question posed by Peter Sarnak, in
[FK] for many quadruples $(k,l_1,l_2,l_3)$
we construct a barrier (with an infinite set $B(\chi)$)
where the imaginary parts of the numbers in the sets $B(\chi)$  are
linearly independent; in particular, we assume all zeros of each
$L(s,\chi)$ are simple, and $L(s,\chi_1)=0=L(s,\chi_2)$ does not occur
for $\chi_1 \ne \chi_2$ and $\Re s >\beta_2$.

So, the results of [FK] show that there are barriers blocking
for large $x$ some ordering of each triple of functions $\pi(x,k,l_i)$
($i=1,2,3$).
 One can ask about barriers blocking some other natural (and usually proven
under ERH$_k$) properties studied in comparative
prime--number theory. We can prove the following results
(a paper is in preparation) concerning the problems 5--7
from section 1:

(1) For many $k,l_1,l_2$, barriers blocking the property
$$
\limsup_{x\to+\infty}\frac{card\{n\le x:\
\pi(n,k,l_1)\ge\pi(n,k,l_2)\}}x>0;
$$

(2) For many $k,l_1,l_2$, barriers blocking the property
$$
\pi(x_\nu,k,l_1)\ge\frac1{\vp(k)}\Li(x_\nu) \text{ and }
\pi(x_\nu,k,l_2)\ge\frac1{\vp(k)}\Li(x_\nu)
$$
for an unbounded sequence $x_\nu$;

(3) For many $k$, barriers blocking the existence
of a sequence $x_\nu\to\infty$ such that
$$
\pi(x_\nu,k,1)<\min_l\pi(x_\nu,k,l) \text{ or }
\pi(x_\nu,k,1)>\max_l\pi(x_\nu,k,l).
$$

%%%%%%%%%%%%%%%%%%%%%%%%%%%%%%%%%%%%%%%%%%%%%%%%%%%%%%%%%%
%
%   REFERENCES
%
%%%%%%%%%%%%%%%%%%%%%%%%%%%%%%%%%%%%%%%%%%%%%%%%%%%%%%%%%%

\enddocument
\bye